\documentclass[12pt]{article}
\usepackage{amsfonts}
\usepackage{amsmath}
\usepackage{graphicx}

\input tcilatex
\input tcilatex
\begin{document}

\author{Ozlem Ersoy Hepson \\
%EndAName
Eski\c{s}ehir Osmangazi University, Faculty of Science and Art,\\
Department of Mathematics-Computer, Eski\c{s}ehir, Turkey}
\title{\textbf{Generation of the Trigonometric Cubic B-Spline Collocation
Solutions for the Kuramoto-Sivashinsky(KS) Equation}}
\maketitle

\begin{abstract}
\noindent A recent type of B-spline functions, namely trigonometric cubic
B-splines, are adapted to the collocation method for the numerical solutions
of the Kuramoto-Sivashinsky equation. Having only first and second order
derivatives of the trigonometric cubic B-splines at the nodes forces us to
convert the Kuramoto-Sivashinsky equation to a coupled system of equations
by reducing the order of the higher order terms. Crank-Nicolson method is
applied for the time integration of the space discretized system resulted by
trigonometric cubic B-spline approach. Some initial boundary value problems
are solved to show the validity of the proposed method.

\end{abstract}

\section{Introduction}

The original form of the Kuramoto-Sivashinsky was constructed to describe
pattern formations and dissipation of them in reaction-diffusion system\cite%
{kuramoto1}. In that study, the reductive perturbation method was
implemented for deriving a scale-invariant part from original macroscopic
motion equations. It was also shown that the Ginzburg-Landau equation can
govern the dynamics near an instability point in many cases. The origin of
persistent wave propagation in reaction-diffusion medium was explored by the
same equation\cite{kuramoto2}. It was also used to explain the chaotic
behavior in a distributed chemical reaction due to the unstable growth of a
spatial inhomogeneity taking place in an oscillating medium\cite{kuramoto3}.
Small model thermal diffusive instabilities in laminar flame fronts can also
be represented by the same equation\cite{sivas1,sivas2}. Nonlinear analysis
of flame front stability assuming stiochiometric composition of the
combustible mixture was also studied with constant-density model of a
premixed flame\cite{sivas3}. The one-dimensional form 
\begin{equation}
u_{t}+uu_{x}+\alpha u_{xx}+\vartheta u_{xxxx}=0  \label{1}
\end{equation}%
of the equation appeared in the study \cite{sivas4}. Hyman and Nicolaenko
characterized the transition to chaos of the solutions by numerical
simulations \cite{hyman1}. The Weiss-Tabor-Carnevale technique applied to
the generalised Kuramoto-Sivashinsky equation to extract some particular
analytical solutions\cite{exa1}. In the related literature, the methods
covering simplest equation, homotopy analysis and $\tanh $ and extended $%
\tanh $ techniques derived to determine solitary wave, or multiple soliton
solutions to the Kuramoto-Sivashinsky equation\cite{exa2,exa3,exa4,exa5}.
\noindent Besides the analytical solutions, many numerical techniques
including Chebyshev spectral collocation\cite{numer1}, finite difference and
collocation \cite{numer2}, quintic B-spline \cite{numer3}, radial basis
meshless method of lines \cite{numer4}, and exponential cubic B-spline
method \cite{numer5} have been applied to derive the numerical solutions to
Kuramoto-Sivashinsky equation.

\noindent \qquad Different from the other B-splines techniques based on
classical polynomial cubic, quartic and quintic B-splines\cite%
{alp1,alp2,alp3} or exponential cubic B-splines \cite{alp4}, the
trigonometric cubic B-spline functions have recently appeared. In this
study, we construct a collocation method based on trigonometric cubic
B-spline functions for some initial boundary value problems for the
Kuramoto-Sivashinsky equation. After reducing the order of the term with the
fourth order derivative to two, we discretize the resultant system by using
Crank-Nicolson method in time. Performing the linearization of the nonlinear
term lead us to discretize the system by trigonometric cubic B-spline
functions. As a result of adapting the initial and boundary conditions, the
iteration algorithm will be ready to run.

To solve the initial value (\ref{1}) numerically we first replace it by a
system which is first order in the time derivative%
\begin{equation}
\begin{array}{r}
u_{t}+uu_{x}+\alpha v+\vartheta v_{xx}=0 \\ 
v-u_{xx}=0%
\end{array}
\label{2}
\end{equation}%
To complete the usual classical mathematical statement of the problem, the
initial and the boundary conditions are chosen as to be%
\begin{equation}
u(x,0)=u_{0}  \label{3}
\end{equation}%
and%
\begin{equation}
\begin{array}{c}
u(x_{0},t)=g_{0},\text{ }u(x_{N},t)=g_{1}, \\ 
u_{x}(x_{0},t)=0,\text{ }u_{x}(x_{N},t)=0, \\ 
u_{xx}(x_{0},t)=0,\text{ }u_{xx}(x_{N},t)=0.%
\end{array}
\label{4}
\end{equation}

\section{Cubic Trigonometric B-spline Collocation Method}

Consider a uniform partition of the problem domain $[a=x_{0},b=x_{N}]$ at
the knots $x_{i},$ $i=0,...,N$ with mesh spacing $h=(b-a)/N.$ On this
partition together with additional knots $x_{-2},x_{-1},x_{N+1},x_{N+2}$
outside the problem domain, $T_{i}(x)$ can be defined as

\begin{equation}
T_{i}(x)=\frac{1}{\gamma }\left \{ 
\begin{tabular}{ll}
$W^{3}(x_{i-2}),$ & $x\in \left[ x_{i-2},x_{i-1}\right] $ \\ 
$%
W(x_{i-2})(W(x_{i-2})Y(x_{i})+Y(x_{i+1})W(x_{i-1}))+Y(x_{i+2})W^{2}(x_{i-1}),
$ & $x\in \left[ x_{i-1},x_{i}\right] $ \\ 
$%
W(x_{i-2})Y^{2}(x_{i+1})+Y(x_{i+2})(W(x_{i-1})Y(x_{i+1})+Y(x_{i+2})W(x_{i})),
$ & $x\in \left[ x_{i},x_{i+1}\right] $ \\ 
$Y^{3}(x_{i+2}),$ & $x\in \left[ x_{i+1},x_{i+2}\right] $ \\ 
$0,$ & $\text{otherwise}$%
\end{tabular}%
\right.   \label{5}
\end{equation}%
where $W(x_{i})=\sin (\frac{x-x_{i}}{2}),\hat{Y}(x_{i})=\sin (\frac{x_{i}-x}{%
2}),\gamma =\sin (\frac{h}{2})\sin (h)\sin (\frac{3h}{2}).$ The twice
continuously differentiable piecewise trigonometric B-spline function set $%
\{T_{i}(x)\}_{i=-1}^{N+1},$ forms a basis for the functions defined in the
same interval \cite{base1,base2}.

$T_{i}(x)$ are twice continuously differentiable piecewise trigonometric
cubic B-spline on the interval $[a,b]$.\ The iterative formula \bigskip

\begin{equation}
T_{i}^{k}(x)=\frac{\sin (\frac{x-x_{i}}{2})}{\sin (\frac{x_{i+k-1}-x_{i}}{2})%
}T_{i}^{k-1}(x)+\frac{\sin (\frac{x_{i+k}-x}{2})}{\sin (\frac{x_{i+k}-x_{i+1}%
}{2})}T_{i+1}^{k-1}(x),\text{ }k=2,3,4,...  \label{6}
\end{equation}%
gives the cubic B-spline trigonometric functions starting with the
CTB-splines of order $1:$

\begin{equation*}
T_{i}^{1}(x)=\left \{ 
\begin{tabular}{c}
$1,\ x\in \lbrack x_{i},x_{i+1})$ \\ 
$0$ \  \ ,otherwise.%
\end{tabular}%
\right.
\end{equation*}%
The graph of the trigonometric cubic B-splines over the interval $[0.1]$ is
depicted in Fig. 1.%
\begin{equation*}
\begin{array}{c}
\FRAME{itbpF}{3.7075in}{3.7075in}{0in}{}{}{fig1.bmp}{\special{language
"Scientific Word";type "GRAPHIC";maintain-aspect-ratio TRUE;display
"USEDEF";valid_file "F";width 3.7075in;height 3.7075in;depth
0in;original-width 3.6599in;original-height 3.6599in;cropleft "0";croptop
"1";cropright "1";cropbottom "0";filename 'Fig1.bmp';file-properties
"XNPEU";}} \\ 
\text{Fig.1: Trigonometric cubic B-splines over the interval }[0,1]%
\end{array}%
\end{equation*}

The nonzero functional and derivative values of trigonometric cubic B-spline
functions at the grids are given in Table 1.

\begin{equation*}
\begin{tabular}{l}
Table 1: Values of $T_{i}(x)$ and its principle two \\ 
derivatives at the knot points \\ 
\begin{tabular}{|l|l|l|l|}
\hline
& $T_{i}(x_{k})$ & $T_{i}^{\prime }(x_{k})$ & $T_{i}^{\prime \prime }(x_{k})$
\\ \hline
$x_{i-2}$ & $0$ & $0$ & $0$ \\ \hline
$x_{i-1}$ & $\sin ^{2}(\frac{h}{2})\csc \left( h\right) \csc (\frac{3h}{2})$
& $\frac{3}{4}\csc (\frac{3h}{2})$ & $\dfrac{3(1+3\cos (h))\csc ^{2}(\frac{h%
}{2})}{16\left[ 2\cos (\frac{h}{2})+\cos (\frac{3h}{2})\right] }$ \\ \hline
$x_{i}$ & $\dfrac{2}{1+2\cos (h)}$ & $0$ & $\dfrac{-3\cot ^{2}(\frac{3h}{2})%
}{2+4\cos (h)}$ \\ \hline
$x_{i+1}$ & $\sin ^{2}(\frac{h}{2})\csc \left( h\right) \csc (\frac{3h}{2})$
& $-\frac{3}{4}\csc (\frac{3h}{2})$ & $\dfrac{3(1+3\cos (h))\csc ^{2}(\frac{h%
}{2})}{16\left[ 2\cos (\frac{h}{2})+\cos (\frac{3h}{2})\right] }$ \\ \hline
$x_{i+2}$ & $0$ & $0$ & $0$ \\ \hline
\end{tabular}%
\end{tabular}%
\end{equation*}%
An approximate solution $U$ and $V$ to the unknown $u$ and $v$ is written in
terms of the expansion of the CTB as

\begin{equation}
U(x,t)=\sum_{i=-1}^{N+1}\delta _{i}T_{i}(x),\text{ }V(x,t)=\sum_{i=-1}^{N+1}%
\phi _{i}T_{i}(x).  \label{7}
\end{equation}%
where $\delta _{i}$ and $\phi _{i}$ are time dependent parameters to be
determined from the collocation points $x_{i},i=0,...,N$ and the boundary
and initial conditions. The nodal values $U$ and its first and second
derivatives at the knots can be found from the (\ref{7}) as 
\begin{equation}
\begin{tabular}{l}
$U_{i}=\alpha _{1}\delta _{i-1}+\alpha _{2}\delta _{i}+\alpha _{1}\delta
_{i+1}$ \\ 
$U_{i}^{\prime }=\beta _{1}\delta _{i-1}-\beta _{1}\delta _{i+1}$ \\ 
$U_{i}^{\prime \prime }=\gamma _{1}\delta _{i-1}+\gamma _{2}\delta
_{i}+\gamma _{1}\delta _{i+1}$%
\end{tabular}%
\begin{tabular}{l}
$V_{i}=\alpha _{1}\phi _{i-1}+\alpha _{2}\phi _{i}+\alpha _{1}\phi _{i+1}$
\\ 
$V_{i}^{\prime }=\beta _{1}\phi _{i-1}-\beta _{1}\phi _{i+1}$ \\ 
$V_{i}^{\prime \prime }=\gamma _{1}\phi _{i-1}+\gamma _{2}\phi _{i}+\gamma
_{1}\phi _{i+1}$%
\end{tabular}
\label{8}
\end{equation}%
\begin{equation}
\begin{array}{lll}
\alpha _{1}=\sin ^{2}(\frac{h}{2})\csc (h)\csc (\frac{3h}{2}) & \alpha _{2}=%
\dfrac{2}{1+2\cos (h)} & \beta _{1}=-\frac{3}{4}\csc (\frac{3h}{2}) \\ 
\gamma _{1}=\dfrac{3((1+3\cos (h))\csc ^{2}(\frac{h}{2}))}{16(2\cos (\frac{h%
}{2})+\cos (\frac{3h}{2}))} & \gamma _{2}=-\dfrac{3\cot ^{2}(\frac{h}{2})}{%
2+4\cos (h)} & 
\end{array}
\label{9}
\end{equation}

When KS equation is space-splitted as (\ref{2}), The system includes the
second-order derivatives so that smooth approximation can constructed with
the combination of the trigonometric cubic B-splines. The time integration
of the space-splitted system (\ref{2}) is performed by the Crank-Nicolson
method as%
\begin{equation}
\begin{array}{r}
\dfrac{U^{n+1}-U^{n}}{\Delta t}+\dfrac{(UU_{x})^{n+1}+(UU_{x})^{n}}{2}%
+\alpha \dfrac{V^{n+1}+V^{n}}{2}+\vartheta \dfrac{V_{xx}^{n+1}+V_{xx}^{n}}{2}%
=0 \\ 
\\ 
\dfrac{V^{n+1}+V^{n}}{2}-\dfrac{U_{xx}^{n+1}+U_{xx}^{n}}{2}=0%
\end{array}
\label{10}
\end{equation}%
where $U^{n+1}=U(x,(n+1)\Delta t)$ represent the solution at the $(n+1)$th
time level. Here $t^{n+1}=t^{n}+\Delta t$, $\Delta t$ is the time step,
superscripts denote $n$ th time level, $t^{n}=n\Delta t.$

One linearize terms $(UU_{x})^{n+1}$and $(UUx)^{n}$ in (\ref{10})\ as \cite%
{rubin}%
\begin{equation*}
\begin{array}{l}
(UUx)^{n+1}=U^{n+1}U_{x}^{n}+U^{n}U_{x}^{n+1}-U^{n}U_{x}^{n} \\ 
(UUx)^{n}=U^{n}U_{x}^{n}%
\end{array}%
\end{equation*}%
to obtain the time-integrated linearized the KS Equation:%
\begin{equation}
\begin{array}{r}
\dfrac{2}{\Delta t}U^{n+1}-\dfrac{2}{\Delta t}%
U^{n}+U^{n+1}U_{x}^{n}+U^{n}U_{x}^{n+1}+\alpha \left( V^{n+1}+V^{n}\right)
+\vartheta (V_{xx}^{n+1}+V_{xx}^{n})=0 \\ 
\\ 
\dfrac{V^{n+1}+V^{n}}{2}-\dfrac{U_{xx}^{n+1}+U_{xx}^{n}}{2}=0%
\end{array}
\label{11}
\end{equation}%
To proceed with space integration of the (\ref{11}), an approximation of $%
U^{n}$ and $V^{n}$ in terms of the unknown element parameters and
trigonometric cubic B-splines separately can be written as (\ref{7}).
Substitute Eqs (\ref{8}) into (\ref{11})\ and collocate the resulting the
equation at the knots $x_{i},$ $i=0,...,N$ yields a linear algebraic system
of equations:

\begin{equation}
\begin{tabular}{l}
\begin{tabular}{l}
$\left[ \left( \frac{2}{\Delta t}+K_{2}\right) \alpha _{1}+K_{1}\beta _{1}%
\right] \delta _{m-1}^{n+1}+\left( \alpha \alpha _{1}+\vartheta \gamma
_{1}\right) \phi _{m-1}^{n+1}+\left[ \left( \frac{2}{\Delta t}+K_{2}\right)
\alpha _{2}\right] \delta _{m}^{n+1}+\left( \alpha \alpha _{2}+\vartheta
\gamma _{2}\right) \phi _{m}^{n+1}$ \\ 
$+\left[ \left( \frac{2}{\Delta t}+K_{2}\right) \alpha _{1}-K_{1}\beta _{1}%
\right] \delta _{m+1}^{n+1}+\left( \alpha \alpha _{1}+\vartheta \gamma
_{1}\right) \phi _{m+1}^{n+1}$ \\ 
$=\frac{2}{\Delta t}\alpha _{1}\delta _{m-1}^{n}-\left( \alpha \alpha
_{1}+\vartheta \gamma _{1}\right) \phi _{m-1}^{n}+\frac{2}{\Delta t}\alpha
_{2}\delta _{m}^{n}-\left( \alpha \alpha _{2}+\vartheta \gamma _{2}\right)
\phi _{m}^{n}+\frac{2}{\Delta t}\alpha _{1}\delta _{m+1}^{n}-\left( \alpha
\alpha _{1}+\upsilon \vartheta _{1}\right) \phi _{m+1}^{n}$%
\end{tabular}
\\ 
\begin{tabular}{l}
$-\gamma _{1}\delta _{m-1}^{n+1}+\alpha _{1}\phi _{m-1}^{n+1}-\gamma
_{2}\delta _{m}^{n+1}+\alpha _{2}\phi _{m}^{n+1}-\gamma _{1}\delta
_{m+1}^{n+1}+\alpha _{1}\phi _{m+1}^{n+1}$ \\ 
$=\gamma _{1}\delta _{m-1}^{n}-\alpha _{1}\phi _{m-1}^{n}+\gamma _{2}\delta
_{m}^{n}-\alpha _{2}\phi _{m}^{n}+\gamma _{1}\delta _{m+1}^{n}-\alpha
_{1}\phi _{m+1}^{n},$ $\  \  \  \  \ m=0...N,$ $n=0,1...,$%
\end{tabular}%
\end{tabular}
\label{12}
\end{equation}%
where%
\begin{equation*}
\begin{array}{l}
K_{1}=\alpha _{1}\delta _{i-1}+\alpha _{2}\delta _{i}+\alpha _{1}\delta
_{i+1} \\ 
K_{2}=\beta _{1}\delta _{i-1}-\beta _{1}\delta _{i+1}.%
\end{array}%
\end{equation*}%
The system (\ref{12}) can be converted the following matrices system;%
\begin{equation}
\mathbf{Ax}^{n+1}=\mathbf{Bx}^{n}  \label{13}
\end{equation}%
where%
\begin{equation*}
\mathbf{A=}%
\begin{bmatrix}
\nu _{m1} & \nu _{m2} & \nu _{m3} & \nu _{m4} & \nu _{m5} & \nu _{m2} &  & 
&  &  \\ 
-\gamma _{1} & \alpha _{1} & -\gamma _{2} & \alpha _{2} & -\gamma _{1} & 
\alpha _{1} &  &  &  &  \\ 
&  & \nu _{m1} & \nu _{m2} & \nu _{m3} & \nu _{m4} & \nu _{m5} & \nu _{m2} & 
&  \\ 
&  & -\gamma _{1} & \alpha _{1} & -\gamma _{2} & \alpha _{2} & -\gamma _{1}
& \alpha _{1} &  &  \\ 
&  &  & \ddots & \ddots & \ddots & \ddots & \ddots & \ddots &  \\ 
&  &  &  & \nu _{m1} & \nu _{m2} & \nu _{m3} & \nu _{m4} & \nu _{m5} & \nu
_{m2} \\ 
&  &  &  & -\gamma _{1} & \alpha _{1} & -\gamma _{2} & \alpha _{2} & -\gamma
_{1} & \alpha _{1}%
\end{bmatrix}%
\end{equation*}

\begin{equation*}
\mathbf{B=}%
\begin{bmatrix}
\nu _{m6} & \nu _{m7} & \nu _{m8} & \nu _{m9} & \nu _{m6} & \nu _{m7} &  & 
&  &  \\ 
\gamma _{1} & -\alpha _{1} & \gamma _{2} & -\alpha _{2} & \gamma _{1} & 
-\alpha _{1} &  &  &  &  \\ 
&  & \nu _{m6} & \nu _{m7} & \nu _{m8} & \nu _{m9} & \nu _{m6} & \nu _{m7} & 
&  \\ 
&  & \gamma _{1} & -\alpha _{1} & \gamma _{2} & -\alpha _{2} & \gamma _{1} & 
-\alpha _{1} &  &  \\ 
&  &  & \ddots & \ddots & \ddots & \ddots & \ddots & \ddots &  \\ 
&  &  &  & \nu _{m6} & \nu _{m7} & \nu _{m8} & \nu _{m9} & \nu _{m6} & \nu
_{m7} \\ 
&  &  &  & \gamma _{1} & -\alpha _{1} & \gamma _{2} & -\alpha _{2} & \gamma
_{1} & -\alpha _{1}%
\end{bmatrix}%
\end{equation*}%
and%
\begin{equation*}
\begin{array}{lll}
\nu _{m1}=\left( \frac{2}{\Delta t}+K_{2}\right) \alpha _{1}+K_{1}\beta _{1}
& \nu _{m4}=\left( \alpha \alpha _{2}+\vartheta \gamma _{2}\right) & \nu
_{m7}=-\left( \alpha \alpha _{1}+\vartheta \gamma _{1}\right) \\ 
\nu _{m2}=\left( \alpha \alpha _{1}+\vartheta \gamma _{1}\right) & \nu
_{m5}=\left( \frac{2}{\Delta t}+K_{2}\right) \alpha _{1}-K_{1}\beta _{1} & 
\nu _{m8}=\frac{2}{\Delta t}\alpha _{2} \\ 
\nu _{m3}=\left( \frac{2}{\Delta t}+K_{2}\right) \alpha _{2} & \nu _{m6}=%
\frac{2}{\Delta t}\alpha _{1} & \nu _{m9}=-\left( \alpha \alpha
_{2}+\vartheta \gamma _{2}\right)%
\end{array}%
\end{equation*}

The system (\ref{13}) consist of $2N+2$ linear equation in $2N+6$ unknown
parameters 
\begin{equation*}
\mathbf{x}^{n+1}=(\delta _{-1}^{n+1},\phi _{-1}^{n+1},\delta _{0}^{n+1},\phi
_{0}^{n+1},\ldots ,\delta _{n+1}^{n+1},\phi _{n+1}^{n+1},).
\end{equation*}
To obtain a unique solution, an additional four constraints\ are needed.
These are obtained from the imposition of the Robin boundary conditions so
that $U_{xx}(a,t)=0,$ $V(a,t)=0$ and $U_{xx}(b,t)=0,$ $V(b,t)=0$ gives the
following equations:%
\begin{equation*}
\begin{array}{l}
\gamma _{1}\delta _{-1}+\gamma _{2}\delta _{0}+\gamma _{1}\delta _{1}=0 \\ 
\alpha _{1}\phi _{-1}+\alpha _{2}\phi _{0}+\alpha _{1}\phi _{1}=0 \\ 
\gamma _{1}\delta _{N-1}+\gamma _{2}\delta _{N}+\gamma _{1}\delta _{N+1}=0
\\ 
\alpha _{1}\phi _{N-1}+\alpha _{2}\phi _{N}+\alpha _{1}\phi _{N+1}=0%
\end{array}%
\end{equation*}

Elimination of the parameters $\delta _{-1},\phi _{-1},\delta _{N+1},\phi
_{N+1},$ from the Eq.(\ref{12}), using the above equations gives a solvable
system of $2N+2$ linear equations including $2N+2$ unknown parameters. After
finding the unknown \ parameters via the application of a variant of Thomas
algorithm, approximate solutions at the knots can be obtained by placing
successive three parameters in the Eq.(\ref{8}).

Initial parameters $\delta _{i}^{0},\phi _{i}^{0},$ $i=-1,\ldots ,N+1$ are
needed to start the iteration procedure (\ref{13}). Thus the following
requirements help to determine initial parameters: 
\begin{equation*}
\begin{array}{l}
U_{xx}(a,0)=0=\gamma _{1}\delta _{-1}^{0}+\gamma _{2}\delta _{0}^{0}+\gamma
_{1}\delta _{1}^{0}, \\ 
U_{xx}(x_{i},0)=\gamma _{1}\delta _{i-1}^{0}+\gamma _{2}\delta
_{i}^{0}+\gamma _{1}\delta _{i+1}^{0}=u_{xx}(x_{i},0),i=1,...,N-1 \\ 
U_{xx}(b,0)=0=\gamma _{1}\delta _{N-1}^{0}+\gamma _{2}\delta _{N}^{0}+\gamma
_{1}\delta _{N+1}^{0}, \\ 
V(a,0)=0=\alpha _{1}\phi _{-1}^{0}+\alpha _{2}\phi _{0}^{0}+\alpha _{1}\phi
_{1}^{0} \\ 
V(x_{i},0)=\alpha _{1}\phi _{i-1}^{0}+\alpha _{2}\phi _{i}^{0}+\alpha
_{1}\phi _{i+1}^{0}=v(x_{i},0),i=1,...,N-1 \\ 
V(a,0)=\alpha _{1}\phi _{N-1}^{0}+\alpha _{2}\phi _{N}^{0}+\alpha _{1}\phi
_{N+1}^{0}%
\end{array}%
\end{equation*}

\section{Numerical tests}

To see versatility of the present method, three numerical examples are
studied in this section. The efficiency and accuracy of the solutions will
be determined by using the global relative error using formula%
\begin{equation}
\text{GRE}=\frac{\dsum \limits_{j=1}^{N}\left \vert
U_{j}^{n}-u_{j}^{n}\right \vert }{\dsum \limits_{j=1}^{N}\left \vert
u_{j}^{n}\right \vert }  \label{GRE}
\end{equation}%
where $U$ denotes numerical solution and $u$ denotes analytical solution.

Numerical solution of KS equation (\ref{1}) is obtained for $\alpha =1$ and $%
\vartheta =1$ with the exact solution given by%
\begin{equation*}
u(x,t)=b+\frac{15}{19}d\left[ e\tanh \left( k\left( x-bt-x_{0}\right)
\right) +f\tanh ^{3}\left( k\left( x-bt-x_{0}\right) \right) \right] 
\end{equation*}%
the initial condition is taken from the exact solution together with
boundary conditions given by (\ref{4}). This example is studied in \cite%
{Xu,numer3,Lai}. The above solution models the shock wave propagation with
the speed $b$ and initial position $x_{0}.$We have considered domain as $%
[x_{0},x_{N}]=[-30,30]$ with time step $\Delta t=0.01$ and number of
partitions as $150$. In order to compare the solutions with \cite{numer3}
and \cite{Lai} we have taken $b=5,$ $k=\frac{1}{2}\sqrt{\frac{11}{19}},$ $%
x_{0}=-12,$ $d=\sqrt{\frac{11}{19}},$ $e=-9,$ $f=11.$ Table 2 gives a
comparison between the global relative error found by our method and by
Quintic B-spline collocation method~\cite{numer3} and by Lattice Boltzmann
method \cite{Lai}.

The numerical results are plotted at different time step for $\Delta t=0.005$
and $N=400$ in Fig. 2 and Fig. 3 shows projection of the solution on the x-t
plane. Solution obtained by trigonometric cubic B-spline collocation method
is very close to the exact solutions due to the global relative error
obtained in Table 2.%
\begin{equation*}
\begin{array}{l}
\text{Table 2: Comparison of global relative error for Example a at
different time }t\text{, }N=150 \\ 
\begin{tabular}{|c|c|c|c|}
\hline
Time($t$) & Present Method & \cite{numer3} & \cite{Lai} \\ \hline
$1$ & \multicolumn{1}{|c|}{$2.98416\times 10^{-5}$} & $3.81725\times 10^{-4}$
& $6.7923\times 10^{-4}$ \\ \hline
$2$ & \multicolumn{1}{|c|}{$7.00758\times 10^{-5}$} & $5.51142\times 10^{-4}$
& $1.1503\times 10^{-3}$ \\ \hline
$3$ & \multicolumn{1}{|c|}{$9.51142\times 10^{-5}$} & $7.03980\times 10^{-4}$
& $1.5941\times 10^{-3}$ \\ \hline
$4$ & \multicolumn{1}{|c|}{$1.79237\times 10^{-4}$} & $8.63662\times 10^{-4}$
& $2.0075\times 10^{-3}$ \\ \hline
\end{tabular}%
\end{array}%
\end{equation*}

\begin{equation*}
\begin{array}{cc}
\FRAME{itbpF}{3.7075in}{3.7075in}{0in}{}{}{fig2.bmp}{\special{language
"Scientific Word";type "GRAPHIC";maintain-aspect-ratio TRUE;display
"USEDEF";valid_file "F";width 3.7075in;height 3.7075in;depth
0in;original-width 3.6599in;original-height 3.6599in;cropleft "0";croptop
"1";cropright "1";cropbottom "0";filename 'Fig2.bmp';file-properties
"XNPEU";}} & \FRAME{itbpF}{2.6757in}{2.7095in}{0in}{}{}{fig3.bmp}{\special%
{language "Scientific Word";type "GRAPHIC";maintain-aspect-ratio
TRUE;display "USEDEF";valid_file "F";width 2.6757in;height 2.7095in;depth
0in;original-width 2.6333in;original-height 2.6671in;cropleft "0";croptop
"1";cropright "1";cropbottom "0";filename 'Fig3.bmp';file-properties
"XNPEU";}} \\ 
\text{Figure 2: Solutions of KS equation} & \text{Figure3: Projected
solutions on }xt-\text{plane}%
\end{array}%
\end{equation*}

\textbf{(b)} This example represents chaotic behaviors with the initial
condition,

\begin{equation*}
u(x,0)=\cos (\frac{x}{2})\sin (\frac{x}{2})
\end{equation*}%
with the boundary condition%
\begin{equation*}
u_{xx}(0,t)=0,\text{ }u_{xx}(4\pi ,t)=0
\end{equation*}%
The computational domain $[x_{0},x_{N}]=[0,4\pi ]$ is used with $N=512,$ $%
\Delta t=0.001,$ $\alpha =1.$ It is shown that KS-Equation is highly
sensitive for choice of the parameter $\vartheta .$ In Figs. 4-7, we can
observe the solution pattern exhibiting complete chaotic behaviors on the $%
xt-$plane, respectively. Figures illustrate that for the smaller value of $%
\vartheta ,$ chaotic behavior starts to evolve earlier and \ seen more
complex instabilities.%
\begin{equation*}
\begin{array}{cc}
\FRAME{itbpF}{4.2376in}{4.1018in}{0in}{}{}{1.jpg}{\special{language
"Scientific Word";type "GRAPHIC";maintain-aspect-ratio TRUE;display
"USEDEF";valid_file "F";width 4.2376in;height 4.1018in;depth
0in;original-width 4.1874in;original-height 4.0525in;cropleft "0";croptop
"1";cropright "1";cropbottom "0";filename '1.jpg';file-properties "XNPEU";}}
& \FRAME{itbpF}{4.2168in}{4.0923in}{0in}{}{}{2.jpg}{\special{language
"Scientific Word";type "GRAPHIC";display "USEDEF";valid_file "F";width
4.2168in;height 4.0923in;depth 0in;original-width 4.1667in;original-height
4.1563in;cropleft "0";croptop "1";cropright "1";cropbottom "0";filename
'2.jpg';file-properties "XNPEU";}} \\ 
\text{Figure 4: Solutions on }xt-\text{plane for }\vartheta =0.05 & \text{%
Figure 5: Solutions on }xt-\text{plane for }\vartheta =0.02%
\end{array}%
\end{equation*}%
\begin{equation*}
\begin{array}{cc}
\FRAME{itbpF}{4.2376in}{4.2168in}{0in}{}{}{3.jpg}{\special{language
"Scientific Word";type "GRAPHIC";maintain-aspect-ratio TRUE;display
"USEDEF";valid_file "F";width 4.2376in;height 4.2168in;depth
0in;original-width 4.1874in;original-height 4.1667in;cropleft "0";croptop
"1";cropright "1";cropbottom "0";filename '3.jpg';file-properties "XNPEU";}}
& \FRAME{itbpF}{4.2168in}{4.2065in}{0in}{}{}{4.jpg}{\special{language
"Scientific Word";type "GRAPHIC";maintain-aspect-ratio TRUE;display
"USEDEF";valid_file "F";width 4.2168in;height 4.2065in;depth
0in;original-width 4.1667in;original-height 4.1563in;cropleft "0";croptop
"1";cropright "1";cropbottom "0";filename '4.jpg';file-properties "XNPEU";}}
\\ 
\text{Figure 6: Solutions on }xt-\text{plane for }\vartheta =0.01 & \text{%
Figure 7: Solutions on }xt-\text{plane for }\vartheta =0.002%
\end{array}%
\end{equation*}

\textbf{(c)} The KS equation (\ref{1}) is obtained for $\alpha =1$ and $%
\vartheta =1$. This example represents the simplest nonlinear partial
differential equation showing chaotic behavior when spatial domain is
finite, with the Gaussian initial condition,%
\begin{equation*}
u(x,0)=-\exp (-x^{2})
\end{equation*}%
with the boundary condition%
\begin{equation*}
u(x_{0},t)=0,\text{ }u(x_{N},t)=0
\end{equation*}%
The computational domain $[x_{0},x_{N}]=[-30,30]$ with $N=120,$ $\Delta
t=0.001.$ In Figs. 8 and 9, we can observe the convergent numerical results
by our trigonometric cubic B-Spline method of lines with complete chaotic
behavior at $t=5$ and $t=20$, respectively. It is observed that the result
shows same characteristics as in \cite{numer3}.%
\begin{equation*}
\begin{tabular}{ll}
\FRAME{itbpF}{3.1652in}{2.6212in}{0in}{}{}{4a.jpg}{\special{language
"Scientific Word";type "GRAPHIC";maintain-aspect-ratio TRUE;display
"USEDEF";valid_file "F";width 3.1652in;height 2.6212in;depth
0in;original-width 5.2295in;original-height 4.3232in;cropleft "0";croptop
"1";cropright "1";cropbottom "0";filename '4a.jpg';file-properties "XNPEU";}}
& \FRAME{itbpF}{3.1592in}{2.6463in}{0in}{}{}{4b.jpg}{\special{language
"Scientific Word";type "GRAPHIC";maintain-aspect-ratio TRUE;display
"USEDEF";valid_file "F";width 3.1592in;height 2.6463in;depth
0in;original-width 5.2191in;original-height 4.3647in;cropleft "0";croptop
"1";cropright "1";cropbottom "0";filename '4b.jpg';file-properties "XNPEU";}}
\\ 
Figure 8: The Chaotic Solution of the KSE $t=5$ & Figure 9: The Chaotic
Solution of the KSE $t=20$%
\end{tabular}%
\end{equation*}


\begin{thebibliography}{99}
\bibitem{kuramoto1} Kuramoto, Y., \& Tsuzuki, T. (1975). On the formation of
dissipative structures in reaction-diffusion systems reductive perturbation
approach. Progress of Theoretical Physics, 54(3), 687-699.

\bibitem{kuramoto2} Kuramoto, Y., \& Tsuzuki, T. (1976). Persistent
propagation of concentration waves in dissipative media far from thermal
equilibrium. Progress of theoretical physics, 55(2), 356-369.

\bibitem{kuramoto3} Kuramoto, Y. (1978). Diffusion-induced chaos in reaction
systems. Progress of Theoretical Physics Supplement, 64, 346-367.

\bibitem{sivas1} Michelson, D. M., \& Sivashinsky, G. I. (1977). Nonlinear
analysis of hydrodynamic instability in laminar flames-II. Numerical
experiments. Acta Astronautica, 4(11-12), 1207-1221.

\bibitem{sivas2} Sivashinsky, G. I. (1977). Nonlinear analysis of
hydrodynamic instability in laminar flames-I. Derivation of basic equations.
Acta astronautica, 4(11-12), 1177-1206.

\bibitem{sivas3} Sivashinsky, G. I. (1980). On flame propagation under
conditions of stoichiometry. SIAM Journal on Applied Mathematics, 39(1),
67-82.

\bibitem{sivas4} Sivashinsky, G. I., \& Michelson, D. M. (1980). On
irregular wavy flow of a liquid film down a vertical plane. Progress of
theoretical physics, 63(6), 2112-2114.

\bibitem{hyman1} Hyman, J. M., \& Nicolaenko, B. (1986). The
Kuramoto-Sivashinsky equation: a bridge between PDE's and dynamical systems.
Physica D: Nonlinear Phenomena, 18(1), 113-126.

\bibitem{exa1} Kudryashov, N. A. (1990). Exact solutions of the generalized
Kuramoto-Sivashinsky equation. Physics Letters A, 147(5-6), 287-291.

\bibitem{exa2} Kudryashov, N. A. (2005). Simplest equation method to look
for exact solutions of nonlinear differential equations. Chaos, Solitons \&
Fractals, 24(5), 1217-1231.

\bibitem{exa3} Abbasbandy, S. (2008). Solitary wave solutions to the
Kuramoto-Sivashinsky equation by means of the homotopy analysis method.
Nonlinear Dynamics, 52(1-2), 35-40.

\bibitem{exa4} Chen, H., \& Zhang, H. (2004). New multiple soliton solutions
to the general Burgers-Fisher equation and the Kuramoto-Sivashinsky
equation. Chaos, Solitons \& Fractals, 19(1), 71-76.

\bibitem{exa5} Wazwaz, A. M. (2006). New solitary wave solutions to the
Kuramoto-Sivashinsky and the Kawahara equations. Applied Mathematics and
Computation, 182(2), 1642-1650.

\bibitem{numer1} Khater, A. H., \& Temsah, R. S. (2008). Numerical solutions
of the generalized Kuramoto-Sivashinsky equation by Chebyshev spectral
collocation methods. Computers \& Mathematics with Applications, 56(6),
1465-1472.

\bibitem{numer2} Lakestani, M., \& Dehghan, M. (2012). Numerical solutions
of the generalized Kuramoto-Sivashinsky equation using B-spline functions.
Applied Mathematical Modelling, 36(2), 605-617.

\bibitem{numer3} Mittal, R. C., \& Arora, G. (2010). Quintic B-spline
collocation method for numerical solution of the Kuramoto-Sivashinsky
equation. Communications in Nonlinear Science and Numerical Simulation,
15(10), 2798-2808.

\bibitem{numer4} Haq, S., Bibi, N., Tirmizi, S. I. A., \& Usman, M. (2010).
Meshless method of lines for the numerical solution of generalized
Kuramoto-Sivashinsky equation. Applied Mathematics and Computation, 217(6),
2404-2413.

\bibitem{numer5} Ersoy, O., \& Dag, I. (2016). The Exponential Cubic
B-Spline Collocation Method for the Kuramoto-Sivashinsky Equation. Filomat,
30(3), 853-861.

\bibitem{alp1} Korkmaz, A., \& Dag, I. (2013). Cubic B-spline differential
quadrature methods and stability for Burgers' equation. Engineering
Computations, 30(3), 320-344.

\bibitem{alp2} Korkmaz, A., \& Dag, I. (2013). Numerical simulations of
boundary-forced RLW equation with cubic b-spline-based differential
quadrature methods. Arabian Journal for Science and Engineering, 38(5),
1151-1160.

\bibitem{alp3} Korkmaz, A., \& Dag, I. (2016). Quartic and quintic B-spline
methods for advection--diffusion equation. Applied Mathematics and
Computation, 274, 208-219.

\bibitem{alp4} Korkmaz, A., \& Akmaz, H. K. (2015). Numerical Simulations
for Transport of Conservative Pollutants. Selcuk Journal of Applied
Mathematics, 16(1).

\bibitem{rubin} S.G. Rubin, R.A. Graves, Cubic spline approximation for
problems in fluid mechanics, Nasa TR R-436, Washington DC, 1975.

\bibitem{base1} Lyche, T., \& Winther, R. (1979). A stable recurrence
relation for trigonometric Bsplines, Journal of Approximation theory, 25(3),
266-279.

\bibitem{base2} Walz, G. (1997). Identities for trigonometric B-splines with
an application to curve design. BIT Numerical Mathematics, 37(1), 189-201

\bibitem{Xu} Y. Xu, C. W. Shu, Local discontinuous Galerkin methods for the
Kuramoto--Sivashinsky equations and the Ito-type coupled KdV equations,
Comput. Meth. Appl. Mech. Eng. 195 (2006) 3430--3447.

\bibitem{Lai} H. Lai, C. Ma, Lattice Boltzmann method for the generalized
Kuramoto--Sivashinsky equation, Physica A 388 (2009) 1405--1412.
\end{thebibliography}
\end{document}